\input amstex
\documentstyle{amsppt}
\magnification1100

\def\bl#1{\smallskip\noindent$\bullet$\quad{\it #1}\quad}
\def\id{\text{id}}
\def\sm{\left(\smallmatrix}
\def\endsm{\endsmallmatrix\right)}

\topmatter
\title Dual 2-complexes in 4-manifolds\endtitle
\author Frank Quinn\endauthor

\date September 2000\enddate
\address 
Mathematics, Virginia Tech, Blacksburg VA 24061-0123\endaddress
\email  quinn\@math.vt.edu\endemail
\abstract This paper concerns decompositions of smooth 4-manifolds as the union of two handlebodies,
each with handles of index~$\leq 2$ (``Heegard'' decompositions). In dimensions $\geq5$ results of
Smale (trivial $\pi_1$) and Wall (general $\pi_1$) describe analogous decompositions up to
diffeomorphism in terms of homotopy type of skeleta or chain complexes. In dimension 4 we show the same
data determines decompositions up to 2-deformation of their spines. In higher dimensions spine
2-deformation
 implies diffeomorphism, but in dimension 4 the fundamental group of the boundary may change. Sample
results: (1.5) Two 2-complexes are (up to 2-deformation) dual spines of a Heegard decomposition of the
4-sphere if and only if they satisfy the conclusions of the Alexander-Lefshetz duality theorem
($H_1K\simeq H^2L$ and
$H_2K\simeq H^1L$).  (3.3) If $(N,\partial N)$ is 1-connected then there is a ``pseudo'' handle
decomposition without 1-handles, in the sense that there is a pseudo collar $(M,\partial N)$ (a
relative 2-handlebody with spine that 2-deforms to $\partial N$) and $N$ is obtained from this by
attaching handles of index $\geq 2$. 
\endabstract
\endtopmatter

\head 1. Results\endhead
By analogy with the 3-dimensional definition, we define a {\it Heegard
decomposition\/} of a 4-manifold to be a description as  a union of two 2-handlebodies along their
boundary. A relative version is defined in \S3. A 2-handlebody has a 2-complex spine. We define a
``tamely embedded'' 2-complex in $N$ to be a spine of a 2-handlebody in $N$, and the handlebody is a
``regular neighborhood'' of the complex. In these terms a Heegard decomposition is a tamely embedded
2-complex whose complement is a regular neighborhood of another tamely embedded 2-complex (the
``dual'').

This paper concerns existence and deformations of
Heegard decompositions, as measured  by chain complexes and 2-complex spines. In these terms they are
fairly flexible: Corollary 1.5, for instance, generalizes the theorem of Huck \cite{H} that any two
acyclic 2-complexes can be  (up to 2-deformation)  dually embedded in $S^4$. Generally they continue
the theme developed by C\. T\. C\. Wall \cite{W1 -- W4} that CW and handlebody structures are
 faithfully described by  cellular chain complexes. This development was motivated by potential
applications to invariants of smooth 4-manifolds, see~\S2.3. 

The paper is organized as follows: Results and
definitions are given in this section. Section 2 lists  facts and questions about Heegard
decompositions. More detailed  relative versions of the main 
results are given in~\S3. Section 4 begins the proof by showing how to align 1-skeleta of chain and CW
complexes. The proof of the realization theorem is given in~\S5.
Section 6  gives a characterization of chain equivalence of 2-complexes in terms of geometric moves.
Finally  the deformation theorem is proved in~\S7.

An 
objective in the following  is to have hypotheses as weak as possible (chain complexes), and
conclusions as strong as possible (2-deformation, which is probably stronger than simple homotopy).
In the statements $N$ is a closed connected smooth 4-manifold, and chain complexes are finitely
generated free based complexes over $Z\pi_1N$. If $K\to N$ is a CW complex then $C^c_*K$ denotes the
cellular chains of the cover of $K$ induced from the universal cover of $N$. Other definitions are
given after the statements, and relative versions are given in~\S3. 

\proclaim{1.1 Theorem (Realization)} Suppose $D_*$ is a chain complex with a chain map $D_*\to
C^c_*(N)$. Then
$D_*$ is simple chain equivalent to the cellular chains of one side of a Heegard decomposition of~$N$ if
and only if it is homologically 2-dimensional and $H_0D\to H_0C^c_*N$ $ (=Z)$ is an
isomorphism.\endproclaim

\proclaim{1.2 Theorem (Deformation)} Suppose $N=M\cup W$ is a Heegard decomposition and $K\to N$ is a
2-complex.
 Then there is an ambient 2-deformation of $W$ to a decomposition $M'\cup W'$ with a
2-deformation
$K\to M'$ if and only if there is a simple chain equivalence $C^c_*K\to C^c_*M$ that chain-homotopy
commutes with the inclusions.
\endproclaim
The balance of this section contains
corollaries and definitions.

\proclaim{1.3 Corollary (Realization of CW spines)}  A finite CW
2-complex $K\to N$ 2-deforms to the spine of half a Heegard decomposition if and only if $K$ is
connected and $\pi_1K\to \pi_1N$ is onto.\endproclaim

\demo{Proof} If $N=M\cup W$ then $W$ is connected and onto $\pi_1N$ because its 1-skeleton is a
1-skeleton of
$N$. For the converse note that connected and onto $\pi_1N$ implies that the cover of $K$ induced from
the universal cover of $N$ is connected. This means $H_0C^c_*K\to H_0C^c_*N\simeq Z$ is an isomorphism.
According to the realization theorem there is a Heegard decomposition $N=M\cup W$ with $C^c_*K\simeq
C^c_*M$ a simple equivalence. But then the deformation theorem shows we can deform $W$ to realize $K$
(up to 2-deformation) as the spine of the complement.\enddemo

In the next corollary we specify the spines of both sides simultaneously. This requires a
 hypothesis that encodes the algebraic duality of the spines in $N$. Suppose $N=M\cup W$ is a Heegard
decomposition. Then there is a diagram of chain complexes,

$$
\CD C_*M@>{\text{dual}}>> C^{4-*}(M,\partial M) @>\text{excision}>>C^{4-*}(N,W)\\
@VVV @. @VVV\\
C_*N   @>\text{dual}>> C^{4-*}N@>=>>C^{4-*}N\\
@. @. @VVV\\
@. @. C^{4-*}W
\endCD$$
in which the rows are chain  equivalences and the right column is homotopy-exact. The
corresponding diagram of homology may be more familiar; there chain equivalence of the rows corresponds
to isomorphism on homology, and homotopy-exactness on the right gives the long exact sequence in
cohomology. 

Consider the
composition from the upper left to the lower right. Homotopy exactness and the fact that it factors
through the upper right complex gives a chain nullhomotopy of this composition. Conversely a chain
nullhomotopy of the composition determines a lift of $C_*^cM@>>> C^{4-*}_cN$ into $C^{4-*}(N,W)$. The
natural nullhomotopy specifies the lift coming from duality and excision, which is a simple
equivalence. 

We abstract this by saying a {\it simple algebraic duality\/} of $K$, $L$ is a chain nullhomotopy of the
corresponding composition, so that the induced lift is a simple equivalence. More explicitly the
nullhomotopy is for the composition
$$C_*^cK@>>> C_*^cN @>\text{dual}>> C^{4-*}_cN @>>> C^{4-*}_c L.$$
The nullhomotopy induces a lift $C_*^cK@>>> C^{4-*}_c(N,L)$, and we require this to be a simple
equivalence. The next corollary shows this characterizes spines of Heegard decompositions.

\proclaim{1.4 Corollary (Characterization of dual spines)} 2-dimensional CW complexes $K, L\to N$
2-deform to spines of dual parts of a Heegard decomposition if and only if there is a  chain
nullhomotopy giving simple algebraic duality on the chain level.\endproclaim

The 2-deformations in the conclusion actually preserve the data in the sense that $K\to M\subset N$ and
$L\to W\subset N$ are homotopic to the original maps, and these homotopies together with the canonical
chain homotopy for $M\cup W$ give (up to chain homotopy) the chain homotopy in the hypothesis. 

\demo{Proof of 1.4} The `only if' part of the statement is explained above. For the converse suppose a
chain homotopy is given. Use the first corollary to realize $K$ as the spine of $M$ in a decomposition
$M\cup W$. The canonical duality for the decomposition gives a simple equivalence $C^c_*K\to
C^{4-*}(N,W)$. The hypothesized duality for $K, L$ gives a simple equivalence $C^c_*K\to C^{4-*}(N,L)$.
Composing one with the inverse of the other and taking duals gives a simple equivalence $C_{*}(N,L)\to
C_{*}(N,W)$. This chain-commutes with the maps from $C_*N$, so via quotients induces a simple
equivalence
$C_{*}L\to C_{*}W$. We can now apply the deformation theorem to ambiently 2-deform $M$ to $M'$ so that
$L$ 2-deforms to the skeleton of the new complement. This gives the desired Heegard decomposition.
\enddemo

The special case of $N=S^4$ can be described explicitly because its homology vanishes, and chain
complexes over the trivial group are determined by homology.

\proclaim{1.5 Corollary (Dual spines in $S^4$)} Two 2-complexes $K$, $L$ occur (up to 2-deformation) as
dual spines of a Heegard decomposition of $S^4$ if and only if they are connected and satisfy
Alexander-Lefshetz duality
$H_1(K)\simeq H^2(L)$ and
$H_2(K)\simeq H^1(L)$.\endproclaim

Huck \cite{H} used a direct explicit construction to realize arbitrary acyclic ($H_1K=H_2K=0$)
2-complexes as spines. His proof could probably be elaborated to prove 1.5. Given a 2-complex $K$ we
note there is a particularly simple choice for dual spine in $S^4$. Suppose
$H_1K\simeq (\oplus_iZ/n_i)\oplus Z^{\beta_1}$ and $H_2(K)\simeq Z^{\beta_2}$. Then the 1-point union
$$L= (\vee_i S^1\cup_{n_i}D^2)\vee (\vee^{\beta_2}S^1)\vee (\vee^{\beta_1}S^2)$$
occurs (after 2-deformation) as a dual spine for $K$.

\demo{Proof of 1.5} Choose (necessarily nullhomotopic) maps $K, L\to S^4$. A standard exact
sequence computation shows there are isomorphisms
$H_*K\simeq H^{4-*}(S^4,L)$.  Chain complexes over
the integers are determined up to  chain homotopy equivalence by  homology, so there is a chain
equivalence $C_*K \simeq C^{4-*}(S^4,L)$. This is simple because the Whitehead group of $Z[1]$ is
trivial. The composition $C_*^cK@>>> C^{4-*}_c L.$ used in 1.4 is nullhomotopic because it factors
through $C_*(S^4)$. As explained before 1.4, a nullhomotopy determines a lift $C_*K \to C^{4-*}(S^4,L)$.
A random homotopy can be changed to one whose lift is the chain equivalence found above. This is a
simple algebraic duality in the sense of 1.4 so implies the existence of a Heegard decomposition with
the desired properties.
\enddemo

 We now provide definitions for terms used in the statements of 1.1 and~1.2.

\subhead 1.6 Definition (homologically 2-dimensional)\endsubhead A finitely generated free chain
complex
$D_*$ over a ring
$R$ is {\it homologically $n$-dimensional\/} if $H_jD=0$ for all $j\geq n+1$ and $H^{j+1}D=0$.

 Wall
\cite{W3, W4} shows that a finitely generated free based complex is simple equivalent to an
$n$-dimensional one if and only if it is homologically $n$-dimensional. He shows further that a finite
CW complex is simple equivalent to an $n$-complex, for $n\neq 2$, if and only if the chain complex of
the universal cover is homologically $n$-dimensional  over the fundamental group. The 2-dimensional
case is still mysterious. Note that 1.1 gives realizations of chain complexes by 2-CW-complexes
(the spine of half the Heegard decomposition), but only by homology equivalence over the fundamental
group of
$N$. These are almost never homotopy equivalences because the fundamental groups are usually different.

\subhead 1.7 Definition (2-deformation)\endsubhead Suppose $K$, $L$ are 2-dimensional CW complexes. A
{\it 2-deformation\/}
$K\to L$ is a sequence of moves starting with $K$ and ending with $L$. Each move is either an
elementary expansion or collapse (of dimension $\leq 2$), or change of attaching map of a 2-cell by
homotopy. In some of the literature these are called ``3-deformations''  because homotopy of
2-cells is viewed as expansion and then collapse of a 3-cell. Our terminology follows Wall \cite{W2}.
Note that a 2-deformation determines a simple homotopy equivalence $K\to L$. The (generalized)
Andrews-Curtis conjecture is that any simple equivalence of 2-complexes comes from a 2-deformation. No
one expects this to be true, however.

In \S6 we define a homological analog, where attaching maps of 2-cells are allowed to change by
homology (in a cover) rather than just homotopy. Unlike geometric 2-deformation there is a global
characterization of the resulting equivalence relation in terms of simple equivalence of chain
complexes.

\subhead 1.8 Definition (ambient 2-deformation)\endsubhead
These are defined for 2-handlebodies embedded in a 4-manifold. They change  spines
by 2-deformation and  preserve Heegard
decompositions in the sense that if the complement of the original is a 2-handlebody then the
complement of the deformation has a corresponding natural 2-handlebody structure. These deformations
usually change the fundamental group of the complement, however. 

An {\it ambient 2-deformation\/} of $M$ in $N$ is a sequence of handle moves in $M$ and homotopically
trivial handle moves over handles in the complement. 
 Handle moves in $M$ include introduction or omission of cancelling of pairs of
handles and do not change the decomposition of $N$. To describe {\it homotopically trivial\/} moves
suppose
 $A$ is a 2-handle in $M$ and $B$ is a
2-handle in the complement attached to $\partial M$ disjointly from $A$. Let $f$, $g$ be embedded paths
from $A$ to
$B$ in $\partial M$ disjoint from
$A\cup B$ except for their endpoints. Suppose also that they are homotopic rel ends by a homotopy in
$M-A$. Obtain a 2-handle $A'$ by pushing $A$ over $B$ twice: once along $f$, and
once along $g$ with opposite sign. Define $M'$ to be $(M-A)\cup A'$. Note that the attaching maps
of $A$ and $A'$ are homotopic in $M-A$ because the curves are homotopic and the signs are opposite.
This gives a 2-deformation between the spines of $M$ and $M'$. The handles of the complement have been
rearranged but no new ones have been added. 

In higher dimensions this move would not change anything
because homotopy in $M-A$ would imply homotopy in $\partial M-A$, and homotopy of paths implies isotopy.
Here however the fundamental groups of the boundary and inside may be different, and the paths are in
the 3-manifold
$\partial M$ so may be knotted.

\head 2. Facts and questions\endhead
In this section we provide background  and pose
further questions
\subhead 2.1 Facts\endsubhead

\bl{Existence}Compact smooth 4-manifolds have Heegard decompositions. A compact smooth manifold $N$ has
a handlebody decomposition that is {\it indexed\/} in the sense that all $k$-handles are attached before
any
$(k+1)$-handles are attached. Take such a decomposition and divide the 2-handles into two sets.
Define $M$ to be the union of the 0- and 1-handles, and the 2-handles in one of these sets. Define $W$
to be the rest of $N$, namely the 4- and 3-handles and the remaining 2-handles. The dual handles give a
description of $W$ as a 2-handlebody, and $N = M\cup W$. 

Note the complete freedom in dividing up the
2-handles between the two sides. The 3-dimensional situation is more determined: since we want each
side to be a 1-handlebody, one side must consist of all of the 0- and 1-handles, while the other
must  be all of the 2- and 3-handles. The extra degrees of freedom in dimension 4 should be thought of
as partial compensation for the greater complexity  of 2-handlebodies. Rather than studying a single
Heegard decomposition, as in dimension 3, we should expect to have to modify it to improve the pieces.

\bl{Smoothness}Donaldson, cf\. \cite{DK}, has shown that not all 4-manifolds have smooth structures,
and when they do the set of such structures is quite mysterious. However (equivalence classes of)
handlebody structures correspond exactly to (diffeomorphism classes of) smooth structures. This is
because handlebodies are determined by the attaching maps in one lower dimension, and attaching maps in
3-manifolds can be uniquely smoothed. Thus a Heegard decomposition encodes a smooth structure, and may
give access to invariants of this structure.

\bl{Pictures}Heegard decompositions  give explicit pictures of smooth 4-manifolds. The
Kirby calculus, extended to include 1-handles, gives explicit pictorial descriptions of 4-dimensional
2-handlebodies in terms of links in the 3-sphere, see
\cite{GS}. The main theorem of the calculus is that the boundaries of two such handlebodies are
diffeomorphic if and only if one link diagram can be deformed to the other using ``Kirby moves.'' The
proof gives a sharper conclusion, namely that a particular sequence of moves determines a particular
diffeomorphism (up to isotopy). A Heegard decomposition, including the glueing map on the boundary, can
therefore be described as two link diagrams together with a sequence of Kirby moves relating them. 

Previous descriptions of 4-manifolds have used link diagrams to represent the whole 2-skeleton. This is
ok when there are no 3-handles, and in principle could be useful even when there are 3-handles because
the attaching maps are fairly rigidly determined. In practice these attaching maps can be very
complicated and hard to see. Whether a Heegard Kirby-move description is more useful remains to be seen.

\subhead 2.2 Cautions\endsubhead

\bl{2-handlebody existence}There are no
  topological criteria for a 4-manifold to be a 2-handlebody. Thus we cannot decompose a manifold
and deduce after the fact that it is a Heegard decomposition. There are easy homological criteria for
$n$-handlebodies when
$n\neq 2$
\cite{W3, W4}, but we expect the analog to be false when $n=2$. In particular there should be smooth
4-manifolds that are simple equivalent to 2-complexes but are not 2-handlebodies.  Note that the simple
equivalence hypothesis does not
 eliminate algebraic topology from the picture since the boundary 3-manifold may have a useful
 fundamental group. Casson has used representations of fundamental groups in compact
connected Lie groups to detect counterexamples to a relative version.

\bl{2-handlebody uniqueness} It seems highly unlikely that 2-handlebody structures are unique in the
sense that one can be deformed to any other through 2-handlebodies, even though the corresponding result
is true for $k$-handlebodies,  $k\neq 2$ \cite{W2}. In other words we expect smooth 4-manifolds with
several different 2-handlebody  ``structures'', so the definition of ``Heegard decomposition'' must
include the handlebody structures, not just the underlying manifolds. 

\bl{2-deformation}The results give embeddings of 2-complexes as spines, but only after 2-deformation.
The original 2-complexes may not be embeddable at all, see \cite{FKT}.  Or they may be embeddable but
with strong restrictions on the possible complementary spines \cite{KT}. In either case preliminary
deformations are necessary to destroy structures that lead to such restrictions.

\bl{Homology equivalence}A homology equivalence of 2-complexes induces a chain homotopy equivalence of
cellular chains, but the converse is not true. In fact a map inducing a homology equivalence must
preserve quite a bit of other structure, such as quotients of fundamental groups by various commutator
subgroups. Therefore ``simple homology equivalence'' cannot be subsituted for ``simple chain
equivalence'' in Theorem~1.2.

\subhead 2.3 Questions\endsubhead

\bl{Reconstruction} In 2.1 we observed that a Heegard decomposition can be pictured as a pair of link
diagrams together with a sequence of Kirby moves going between them. We need an explicit equivalence
relation (in terms of ``moves'') on such data so that equivalence classes of data correspond to
diffeomorphism classes of manifolds. Two sorts of moves are needed: one would show how to transfer a
2-handle from one  side of the decomposition to the other. The other sort should explain how to relate
two sequences of Kirby moves that yield isotopic diffeomorphisms of the boundary 3-manifold. The second
part may be more accessible than one might expect.  The {\it
existence\/} of a sequence of Kirby moves corresponds to the existence of a diffeomorphism between
3-manifolds, and this seems to be algorithmically unsolvable in terms of the data given. However {\it
equivalence\/} of two sequences would correspond to automorphisms of a fixed 3-manifold. These are
frequently quite limited and may be detectable. 

\bl{Geometric structures} Eliashberg has observed that a Stein surface has a 2-handlebody structure
and Gompf has characterized the handlebodies obtained this way, see \cite{GS}.  Do all smooth
4-manifolds have Heegard decompositions with Stein pieces? More precisely (in the spirit of 1.2) how
much does one side of a decomposition have to be changed to make the other side Stein? The boundary of
a Stein surface has a contact structure. Are there  decompositions with Stein pieces in which the
glueing map has some nice relationship with the contact structures?  The goal is to find a
comprehensible link between a combinatorial structure and Donaldson or Seiberg-Witten invariants.  The
combinatorial structure would be a Heegard decomposition with Stein pieces, together with a sequence of
restricted Kirby moves giving a ``contact-nice'' glueing map. The link would go through gauge theory on
Stein surfaces. 

\bl{Higher-order deformations} The refined duality results of Krushkal and Teichner \cite{KT} suggest
there should be a similarly refined version of the results of this paper. The analog of \S6 would relate
refined chain equivalence and geometric moves. If $K$ is a 2-complex over a group $\pi$ the algebraic
version  presumably would concern equivalence over $\pi_1K$ modulo the the $k^{\text{th}}$ commutator of
the kernel of $\pi_1K\to \pi$. The homotopy moves would use gropes of height $k$ mapping into the
1-skeleton, with $\pi_1$ going trivially into $\pi$. Finally ``ambient deformation'' would involve
additions along curves that bound a $k$-stage grope with embedded body. 

Presumably something like Milnor invariants would appear as obstructions to $k$-equivalent decompositions
being $k+1$-equivalent.

\head 3. Relative versions\endhead

A  Heegard decomposition of a compact 4-manifold {\it with boundary\/} is a
decomposition $(N,\partial N)=(M, A)\cup (W, V)$ into submanifolds interesecting in a submanifold of
their boundaries, and $(M,A)$, $(W,V)$ are given as relative
2-handlebodies. A relative 2-handlebody structure on
$(M,A)$ is a description as a collar $A\times I$ and handles of index $\leq 2$ attached on $A\times
\{1\}$ and handles of lower index. The result has $A\times \{0\}$ embedded in its boundary, and this is
identified with $A\subset \partial M$.

In the following $N$ is a compact smooth 4-manifold and ``chain complex'' means finitely generated free
based complex over
$Z[\pi_1N]$. The absolute versions focus on geometric conclusions from algebraic hypotheses.
Appropriate formulations of the algebraic hypotheses are necessary as well as sufficient, but trying
to include this in the statements makes them unnecessarily complicated.  

\proclaim{3.1 Theorem (Relative realization)} Suppose a decomposition into codimension 0 submanifolds
$A\cup V=\partial N$ is given,
$D_*$ is a homologically 2-dimensional chain complex, and $D_*\to C^c_*(N,A)$  induces isomorphism on
$H_0$, epimorphism on $H_1$. Then there exists a Heegard decomposition
$(N,\partial N)=(M,A)\cup (W,V)$ with a simple chain equivalence $D_*\to C^c_*(M,A)$\endproclaim
\proclaim{Refinement}  The composition $D_*\to C^c_*(M,A)\to C^c_*(N,A)$ is 
chain homotopic to the original chain
map.\endproclaim

The point of the refinement is that the geometric decomposition realizes the data:
we only have to change the chain map by homotopy to geometrically realize it. It  sometimes happens
in low dimensions that existence of data implies geometric conclusions, but corresponding to different
data. The discrepancy between given and realized data then gives rise to higher-order obstructions.
Fortunately this does not happen here.

Setting $\partial N=\emptyset$ in 3.1 gives~1.1, except for the $H_1$ hypothesis. However the homology
is with $Z\pi_1N$ coefficients, or in other words homology of the universal cover of $N$. This is
simply connected so $H_1$ vanishes, and the hypothesis holds trivially.

To prepare for the deformation statement we recall that ambient deformation does not change
the chains of the complement, in the sense that if $(W,B)$ ambiently 2-deforms to $(W'B)$, and the
complements are $(M,A)$ and $(M',A)$ respectively, then there is a canonical (up to chain homotopy)
simple chain equivalence $C^c_*(M,A)\to C^c_*(M',A)$, together with a chain homotopy making the
diagram induced by inclusions 
$$\CD C^c_*(M,A)@>>> C^c_*(M',A)\\
@VVV@VVV\\
C^c_*(N,A)@>=>> C^c_*(N,A)\endCD$$
commute. This can be seen either directly from the definition of ambient 2-deformation or by using
duality and the corresponding fact for the 2-deformation $W\to W'$.

\proclaim{3.2 Theorem (Relative deformation)} Suppose $(N,\partial N)=(M,A)\cup (W,V)$ is a Heegard
decompositon,
$(K,A)\to (N,A)$ is a map of a  relative 2-complex, and there is a simple chain equivalence
$C^c_*(K,A)\to C^c_*(M,A)$ together with a chain homotopy making the diagram induced by inclusions
commute. Then there is an ambient 2-deformation of $W$ to give $(M',A)\cup(W',V)$, and a relative
2-deformation from $(K,A)$ to the spine of $(M',A)$.
\endproclaim
\proclaim{Refinement} There are
\roster\item a homotopy from the composition of the 2-deformation with the inclusion,
$(K,A)\to(M',A)\to (N,A)$, to the original map;
\item a chain homotopy from  given chain equivalence $C^c_*(K,A)\to C^c_*(M,A)$ to the
the composition  $C^c_*(K,A)\to C^c_*(M',A)\to C^c_*(M,A)$ , where the first map is induced by the
2-deformation to the skeleton and the second is dual to the ambient 2-deformation of $W$; and
\item note that combining the homotopies of (1) and (2) give a chain homotopy making the diagram induced
by inclusions commute. Then there is a chain homotopy between this homotopy and the one given in the
data.
\endroster
\endproclaim
Again the refinement asserts that the output realizes the input data, so no higher obstructions
arise from discrepancies between the two. From a technical point of view the extra precision becomes
important in inductions or other arguments applying the theorem several times.

To illustrate use of the relative versions  we give a version of ``geometrical connectivity'' for
4-manifolds. A special case of Wall's results
\cite{W1} is that in higher dimensions a 1-connected manifold pair is obtained from a collar of the
submanifold of the boundary by adding handles  of index $\geq 2$.  Casson has shown that this is
not generally true for 4-manifolds. However there is a version where the handles  are
added to something which is not quite a collar on the boundary. 
\proclaim{3.2 Corollary} Suppose $N$ is a compact smooth 4-manifold, $A\subset\partial N$ is a
 codimension 0 submanifold, and $(M,A)$ is
relatively 1-connected. Then there is a ``fake collar'', ie a codimension 0 submanifold $M\subset N$
with $M\cap \partial N=A$ and $(M,A)$ a relative 2-handlebody whose skeleton 2-deforms rel $A$ to $A$,
so that $N$ is obtained from $M$ by adding handles of index $\geq 2$. 
\endproclaim
\demo{Proof} Let $V$ denote the closure of the complement of $A$ in $\partial N$. The 1-connected
hypothesis is equivalent to
$C_*(N,V)$ being homologically 2-dimensional. Applying the realization theorem gives a relative Heegard
decomposition $N=M\cup W$ with the inclusion $C_*(M,V)\to C_*(N,V)$ a simple chain equivalence over
$Z\pi_1N$. The inclusion of the trivial relative 2-complex $(A,A)\to (M,A)$ gives a simple equivalence
of chains, so by the deformation theorem $W$ can be ambiently 2-deformed to give $M$ whose skeleton
relatively 2-deforms to $A$. The duals of the handles in $(W,V)$ describe $N$ as obtained from $M$ by
adding handles of index $\geq 2$, so this decomposition satisfies the conclusions of the
Corollary.\enddemo

Note that  $(W,V)\to (N,V)$ is a homology equivalence over $Z\pi_1N$. This is a homotopy equivalence
if and only if
$\pi_1W\to
\pi_1N$ is an isomorphism. In this case $(N,W)$ has the homotopy type of a 2-complex, as well as being
homologically 2-dimensional. Conversely if
$(N,V)$ were known to have the homotopy type of a relative 2-complex then we could realize that complex
(up to 2-deformation) as the spine of $W$, and get $(W,V)\to (N,V)$ a homotopy equivalence.

We close the section with a few observations about the fake collars appearing in the corollary.
First, 2-deformations can be realized by handle moves in a 5-manifold \cite{AC}, so $(M\times I, 
A\times I)\simeq (A\times I^2,A\times I\times\{0\})$. The second observation is that $(M,A)$ is a simple
homology H-cobordism with coefficients $\pi_1M=\pi_1A$, from $A$ to the closure of its complement in
$\partial M$. This is usually not an s-cobordism because $\partial M-A$ usually has different
fundamental group. This must happen if $\pi_1W\to \pi_1N$ is not an isomorphism, and in particular if
$(W,V)$ is not equivalent to a 2-complex.

\head 4. 1-skeleton alignment\endhead
The objective here is to modify chain complexes to standardize their 1-skeletons. If the complex is
cellular chains of a CW complex we want to achieve this standardization through 2-deformations of the
CW complex. The result is similar to the low-dimensional cases of Wall \cite{W2}, see particularly Lemma
3B. 

Fix a group $\pi$. ``Chain complex'' will mean a finitely generated free based $Z[\pi]$ complex which
vanishes in negative degrees. CW complexes will be assumed to come with a homomorphism $\pi_1K\to \pi$,
and maps, deformations, etc\. of CW complexes are understood to commute with these homomorphisms. The
``cellular chains'' of a CW complex are the cellular chains of the induced cover with covering group
$\pi$, with a free basis obtained by lifting cells. 

\proclaim{Lemma (1-skeleton alignment)} Suppose $f\:C_*\to D_*$ is a chain map of free based
$Z[\pi]$ complexes which is an isomorphism on $H_0$ and an epimorphism on $H_1$. Then there is a chain
homotopy commutative diagram 
$$\CD C_*@>f>> D_*\\
@VVV @VVV\\
C'_* @>{f'}>>D'_*
\endCD$$
so the vertical maps are simple equivalences and isomorphisms in degrees $\geq 3$, and $f'$ is a
basis-preserving isomorphism in degrees 0, 1. 

If $D_*$ is the cellular chains of a  CW complex or pair (resp. 4-d handlebody) with connected $\pi$
cover then
$D_*\to D'_*$ can be arranged to be the map on cellular chains induced by a 2-deformation (resp.
handlebody moves). 

If both $C_*$ and $D_*$ are cellular chains of CW complexes or pairs (resp. 4-d handlebodies) with
connected
$\pi$ covers and the isomorphism on $H_0$ is the identity, then both
$C_*\to C'_*$ and $D_*\to D'_*$ can be arranged to be induced by 2-deformation (resp. handle moves).
\endproclaim
In the last case recall that $H_0$ of a connected CW complex is canonically identified with $Z$. A
chain map of cellular chains that induces $-1$ on $H_0$ with respect to these identifications cannot
be made basis-preserving in degree 0 because there are no orientation-reversing endomorphisms of a
point. In the relative case $H_0=0$ so the hypothesis is automatic.

\subsubhead{Proof}\endsubsubhead We prove this by collecting data and then explicitly giving chain
complexes, homotopies, etc. The first step is to align the 0-skeletons, ie\. arrange $C_0\to D_0$ to be
an isomorphism. Assuming this we then show how to align the 1-skeleton. 

There are several procedures for aligning skeleta, depending on dimension and whether we want to realize
the changes by 2-deformation. The simplest procedure fixes a dimension and changes both complexes one 
dimension higher. The bound on the changes means it can be used for 1-skeleta. This does not
work for 0-skeleta when one of the complexes is geometric because the changes cannot be realized
geometrically. In that case we give a more elaborate argument that leaves one complex (the geometric one)
unchanged and modifies the other. This argument  changes the complex two dimensions above the target
dimension, so cannot be used for 1-skeleta. 

\subsubhead{Semi-algebraic 0-skeleton alignment}\endsubsubhead Here we are not trying to
realize  changes in the domain 
 by 2-deformation. The range complex is unchanged, so can be geometric. The data is a chain map of
complexes
$$\CD C_2@>f_2>> D_2\\
@VV{\partial}V @VV{\partial}V\\
C_1@>f_1>> D_1\\
@VV{\partial}V @VV{\partial}V\\
C_0@>f_0>> D_0
\endCD$$
Since this is onto $H_0$ the map $\sm f_0&\partial\endsm\: D_0\oplus
C_1\to D_0$ is onto. Since $D_0$ is free there is  a splitting $\sm g_0\\
s\endsm$. In other words 
$$f_0g_0+\partial\,s =\id.\tag{1}$$ 
A similar argument using the fact that $f_*$ induces an isomorphism on $H_0$ shows there is $t\:C_0\to
C_1$ satisfying
$$g_0f_0 +\partial\,t =\id.\tag{2}$$
In the next step we improve this $t$. Consider 
$$\CD @. C_1\oplus D_2\\
@. @VV{\sm f_1& \partial\endsm}V\\
C_0@>s f_0-f_1t>> D_1
\endCD$$
We claim there is a lifting of the horizontal map to $\sm \alpha\\u\endsm\:C_0\to
(\text{ker}\,\partial)\oplus D_2$. It suffices to show the image of the horizontal is in the image of
this subspace. Since 
$f_*$ is onto on
$H_1$, the image of
$(\text{ker}\,\partial)\oplus D_2\to D_1$ is the kernel of $\partial\:D_1\to D_0$. Therefore to verify
the image of the horizontal map is in this image it suffices to verify the compostion with $\partial$
vanishes. But $\partial(s f_0-f_1t)=(\partial s)f_0 -f_0(\partial t)$, which vanishes by (1) and (2). 

Replace $t$ by $t+\alpha$. Since $\alpha$ has image in $\text{ker}\,\partial$,  property (2) still
holds. The factorization property becomes 
$$\partial\, u+f_1t=s\,f_0.\tag{3}$$
Now form the diagram
$$\CD C_2\oplus C_0 @>{\sm f_2&u\endsm} >> D_2\\
@VV{\sm \partial& -t\\ 0 & f_0\endsm}V @VV{\partial}V\\
C_1\oplus D_0 @>{\sm f_1&s\endsm} >> D_1\\
@VV{\sm f_0\partial &\partial s\endsm}V @VV{\partial}V\\
D_0 @>{\id}>> D_0
\endCD$$
The vertical sequences continue up with higher-degree terms of $C_*$, $D_*$, if there are any.
Identity (3) shows the left vertical composition is trivial, so gives a chain complex, and also that the
upper square commutes so the horizontal maps give a chain map. 

Define the left-hand complex to be $C'_*$ and
the map to be $f'_*$. Define $D'_*=D_*$; since it is unchanged the change is realized by a
2-deformation when $D_*$ is geometric. $f'_*$ is a based isomorphism in degree 0. To finish the
construction we need a simple equivalence $C_*\to C'_*$.

There is a natural chain map of $C_*$ into $C'*$,
$$\CD C_2 @>{\sm \id\\0\endsm}>> C_2\oplus C_0\\
@VV{\partial}V @VV{\sm \partial& -t\\ 0 & f_0\endsm}V\\
C_2 @>{\sm \id\\0\endsm}>>C_1\oplus D_0\\
@VV{\partial}V @VV{\sm f_0\partial &\partial s\endsm}V\\
C_0 @>f_0>> D_0
\endCD$$
It is easy to see that composition of this with $f'_*$ gives $f_*$, so it suffices to show
this is a simple equivalence. This can be done by adding the trivial complex $\id\:C_0\to C_0$ in
degrees 0 and 1, then doing elementary moves to transform the result to $C_*$ plus the trivial complex
$$C_0@>{\sm \id\\0\endsm}>>C_0\oplus D_0 @>{\sm 0&\id\endsm}>>D_0.$$
We omit the details of this.

\subsubhead{Geometric 0-skeleton alignment}\endsubsubhead
Here both $C_*$ and $D_*$ are geometric, and we want the modifications realized by geometric
2-deformations. We begin with the non-relative case. Since both complexes are connected there are
2-deformations to CW complexes with single 0-cells. In the manifold case there are handle moves to
handlebody structures with single 0-handles. Let
$C'_*$, $D'_*$ be the chains of these complexes. 

Since the $\pi$ covers of the spaces are connected the single basepoints give a single copy of $Z\pi$
in degree 0. Thus the chain map has the form 
$$\CD C_1@>f_1>> D_1\\
@VV{\partial}V @VV{\partial}V \\
Z\pi@>f_0>> Z\pi.
\endCD$$
By hypothesis $f_*$ induces the identity on $H_0$. So does the identity map on $Z\pi$, so the image of
$(\id-f_0)\:Z\pi\to Z\pi$ lies in the image of the boundary in $D_*$. Let $s\:Z\pi\to C_1$ be a lift.
Then $s$ defines a chain homotopy of $f_*$ to the chain map 
$$\CD C_1@>f_1+s\partial>> D_1\\
@VV{\partial}V @VV{\partial}V \\
Z\pi@>\id>> Z\pi.
\endCD$$
This is a basis-preserving isomorphism in degree 0 so serves as $f'_*$.

We now consider the relative case, where $C$ and $D$ are the relative chains of pairs. Since the spaces
are connected there are 2-deformations to relative 2-complexes that have no 0-cells. Since the chain
groups in dimension 0 are now trivial the chain map trivially induces an isomorphism. However there is a
little more to do for pairs $(K,A)$ when $A$ is not connected. The condition actually used in the
1-skeleton argument is that all 1-cells are attached at a single point. If $A$ is not connected choose a
linear order of the components of $A$ and deform the CW structure so there is a 1-cell joining each
component to the next one. All other 1-cells can be assumed to be attached at a point.  The  alignment
 needed  is to deform the other CW complex to have a similar structure, so that the chain map preserves
the classes of the arcs connecting pieces of $A$. This begins by finding maps of arcs rel endpoints into
the complex rel $A$ that are liftings on the chain level. The mapping cylinder of these maps has
a CW structure with new 1- and 2-cells. This 2-deforms to the complex, but now has 1-cells in the right
place algebraically. Next we  manipulate the CW structure to get all other 1-cells to be attached at
a point. This can be done to end up with the desired alignment of the special 1-cells. We omit details
since the argument is basically the same as others in this section. 

\subsubhead{Algebraic 1-skeleton alignment}\endsubsubhead The hypotheses are as in the lemma, and
additionally $f_0$ is a basis-preserving isomorphism. To reflect this we denote the bottom term of
$D_*$ also by $C_0$, and $f_0=\id$. Here we do algebraic moves to get a basis-preserving isomorphism in
degree 1 as well.  First we collect some data.

Since $f_0$ and $H_0(f_*)$ are isomorphisms, $f_0$ induces an isomorphism on the images $\partial
C_1\to \partial C_2$. This together with the hypothesis that $H_1(f_*)$ is onto shows that $f_1\:C_1\to
D_1/\partial D_2$ is onto. This means there is a lift $g\:D_1\to C_1$ so $f_1g=\id$ modulo $\partial
D_2$. This in  turn implies a lift $s\:D_1\to D_2$ so 
$$\partial s = \id -f_1g.\tag{1}$$
Now consider the diagram
$$\CD C_2\oplus D_1 @>{\sm 1&0\\0&1\endsm}>> C_2\oplus D_1@>{\sm \partial&-g\\f&s\endsm}>>C_1\oplus D_2
@<{\sm 1&0\\0&1\endsm}<<C_1\oplus D_2\\
@VV{\sm \partial&0\\0&1\endsm}V@VV{\sm \partial&-g\\0&1\endsm}V@VV{\sm \partial&0\\-f&1\endsm}V@VV{\sm
\partial&0\\0&1\endsm}V\\
C_1\oplus D_1 @>{\sm 1&-g\\0&1\endsm}>> C_1\oplus D_1@>{\sm 1&0\\0&1\endsm}>>C_1\oplus D_1
@<{\sm 1&0\\-f&1\endsm}<<C_1\oplus D_1\\
@VV{\sm \partial&0\endsm}V@VV{\sm \partial&\partial\endsm}V@VV{\sm \partial&\partial\endsm}V@VV{\sm
0&\partial\endsm}V\\
C_0@>>>C_0@>>>C_0@<<<C_0
\endCD$$
The identity (1) shows that the columns are all chain complexes and the horizontal maps give chain
maps. Let the second and third columns be $C'_*$ and $D'_*$ respectively, and the map between them
$f'_*$. Note $f'_*$ is a basis-preserving isomorphism in degrees 0 and 1. The first column is a
stabilization of
$C_*$ and the map to the second column is a simple equivalence. Thus this gives the simple equivalence
$C_*\to C'_*$. Similarly the right column is a stabilization of $D_*$ and the map to the third column is
simple. This gives
$D_*\to D'_*$.  It only remains to check that the square obtained by adding the original $f_*\:C_*\to
D_*$ chain-homotopy commutes. This too is a simple consequence of (1). These data therefore give the
conclusion of the lemma.

\subsubhead{Geometric 1-skeleton alignment}\endsubsubhead
By construction the chain maps $C_*\to C'_*$ and $D_*\to D'_*$ consist of stabilizations adding
cancelling 1- and 2-dimensional generators, and endomorphisms using elementary (upper or lower
triangular) matrices on the 1-dimensional generators. Stabilization can be realized by elementary
expansions, or introduction of cancelling 1, 2-handle pairs in the handlebody case. 1-cells
or 1-handles all attached to a single 0-cell or 0-handle can be moved over each other by homotopy (resp.
isotopy) to realize endomorphisms the chain complex by elementary matrices. Therefore if one or both
complexes starts out as cellular chains then the modified complexes are realized by 2-deformation
(resp. handle moves). 

This completes the proof of the 1-skeleton alignment lemma.

\head 5. Proof of the realization theorem\endhead
Suppose, as in the statement of 3.1, that $D_*$ is a homologically 2-dimensional $Z\pi_1N$-complex and
$D_*\to C^c_*(N,A)$ induces isomorphism on $H_0$, epimorphism on $H_1$. According to Wall \cite{W3,
W4}, it is simple equivalent to a 2-dimensional chain complex. Use the 1-skeleton alignment lemma to
change
$D$ by simple equivalence to another 2-dimensional complex, and the handle decomposition on $(N,A)$ by
handle moves, so $D_*\to C^c_*(N,A)$ is a basis-preserving isomorphism in degrees 0 and 1.  

We claim the handlebody structure can be further changed so $D_2\to C^c_2(N,A)$ is a basis-preserving
isomorphism on a summand. Assuming this we complete the proof of the theorem. Simply take $M$ to be the
union of the 0-handles, the 1-handles, and the 2-handles whose corresponding basis elements in
$C^c_2(N,A)$ come from $D_2$. Then $C^c_*(M,A)$ is the image of $D_*$ in $C^c_*(N,A)$, and the map
$D_*\to C^c_*(M,A)$ is a basis-preserving isomorphism and therefore a simple equivalence. Since $M$
contains all the 0- and 1-handles the complement is a 2-handlebody, and gives a Heegard decomposition.

We now prove the claim. The situation is  a chain map which is a based isomorphism in degrees 0, 1. 
$$\CD @. D_2@>>>D_1@>>>D_0\\
@. @VV{f}V@VV=V@VV=V\\
C_3(N,A)@>>>C_2(N,A)@>>>C_1(N,A)@>>>C_0(N,A)
\endCD$$
Introduce in $N$ cancelling 2, 3-handle pairs indexed by the basis of $D_2$. Do elementary moves among
the 2-handles to change the cellular chains by the chain map
$$\CD C_3\oplus D_2@>>{\sm\partial &0\\0&1\endsm}> C_2\oplus D_2 @>>{\sm \partial&0\endsm}>
C_1@>>>\cdots\\ @VV{\sm 1&0\\0&1\endsm}V@VV{\sm 1&-f\\0&1\endsm}V@VV=V\\
 C_3\oplus D_2@>>{\sm\partial &-f\\0&1\endsm}>
C_2\oplus D_2 @>>{\sm \partial&\partial\endsm}> C_1@>>>\cdots
\endCD$$
The composition of the original map from $D$ into this gives
$$\CD @. D_2@>>>D_1@>>>\cdots\\
@. @VV{\sm f\\0\endsm}V@VV=V\\
C_3\oplus D_2@>>{\sm\partial &-f\\0&1\endsm}>
C_2\oplus D_2 @>>{\sm \partial&\partial\endsm}> C_1@>>>\cdots
\endCD$$
But now $\sm 0\\1\endsm\:D_2\to C_3\oplus D_2$ is a chain homotopy that changes this chain map to be
the inclusion of the second summand in dimension 2. This is a basis-preserving injection to a based
summand of the 2-chains, proving the claim and the theorem.

\head 6. Homological 2-deformation\endhead
This is a characterization of simple chain equivalence of 2-complexes in terms of geometric moves. To
some extent it is a subsitute for the unknown and probably false characterization of simple homotopy by
2-deformation. This provides a new derivation and wider context for the ``s-moves'' shown to give simple
homotopy in~\cite{Q}. Homological 2-deformations are used in  the proof of the deformation theorem~3.2.

\subhead Definition (Homological 2-deformation)\endsubhead Suppose $\pi$ is a group and $K$ is a
2-complex with a homomorphism $\pi_1K \to \pi$. A homological 2-deformation (with $Z\pi$ coefficients)
is a sequence of moves, each of which is either a 2-deformation or change of attaching map of a 2-cell
by $Z\pi$ homology in the 1-skeleton.

We describe this last move more explicitly. Suppose $\alpha, \beta\:S^1\to K^1$ are maps with the same
basepoint. Here $K^1$ denotes the 1-skeleton of $K$. These maps differ by $Z\pi$ homology in the
1-skeleton if the composition
$\alpha^{-1}\beta$ extends to a map of an orientable surface with one boundary component $S\to K^1$, and
this map lifts to  the $\pi$ cover of $K^1$. The lifting property is equivalent to triviality of the
composition on $\pi_1$,
$\pi_1S\to \pi_1K^1\to \pi$.

\proclaim{Proposition (Homological deformation is chain equivalence)} Suppose $\pi$ is a group and $K$,
$L$ are 2-complexes with homomorphisms $\pi_1\to \pi$. A homological 2-deformation $K\to L$ over $\pi$
induces a simple chain equivalence of $Z\pi$ cellular chain complexes $C^c_*K\to C^c_*L$. Conversely if
$K$ and $L$ are connected and the homomorphisms to $\pi$ are onto then any simple chain equivalence is
realized (up to chain homotopy) by a homological 2-deformation.
\endproclaim
A straightforward relative version also holds.

\demo{Proof}
For the first statement recall that 2-deformations induce simple chain
equivalences over $\pi_1$ and therefore over any other coefficients. Changing attaching map of a 2-cell
by homology does not change the cellular chain complex at all. The only issue is the boundary
homomorphism $C_2\to C_1$, which is defined using homology classes in the 1-skeleton. Therefore a
homological 2-deformation also induces a simple chain equivalence.

For the other direction first apply the 1-skeleton alignment lemma. This gives 2-deformations to $K'$
and $L'$ together with a chain map of cellular chain complexes (over $\pi$) that is a basis-preserving
isomorphism on $C_0$ and $C_1$. Since it is a chain homotopy equivalence and the complexes are
2-dimensional it must be an isomorphism on $C_2$ as well. Since the equivalence is simple this
isomorphism is also simple. This means that (possibly after stabilization) it is a product of elementary
matrices and diagonal matrices with group elements on the diagonal. Such a matrix can be realized by
homotopies of attaching maps and changing choice of liftings used as bases. Doing these moves to $K$
gives $K$, $L$ and a basis-preserving isomorphism of cellular chains.

There is a unique (up to homotopy) way to identify the 1-skeletons of $K$, $L$ so that the identification
induces the chain isomorphism. We now regard $K$, $L$ as 2-complexes with the same 1-skeleton. The chain
isomorphism gives a bijection between 2-cells so corresponding cells have the same boundaries in
$C^c_1K^1$. But $C^c_1$ is defined to be $H_1(K^1,K^0;Z\pi)$, so the attaching maps are homologous in the
$\pi$ cover of $K^1$. Equivalently the difference is nullhomologous. But a map of $S^1$ into a space is
nullhomologous if and only if it extends to a map of an oriented surface. Therefore the attaching maps in
the two complexes differ by oriented surfaces as specified in the definition, so they are obtained from
each other by homological 2-deformation. 

Note the final homological moves realize the chain map (which is by then the identity). This chain map is
homotopic to the original, so the original is realized up to homotopy.
\enddemo

\subhead 6.2 Simple homotopy equivalence\endsubhead A simple {\it homotopy\/} equivalence of {2-}
complexes
$K\to L$ is a map which is an isomorphism on $\pi_1$ and a simple equivalence on cellular chain complexes
over
$\pi_1L$. This is in particular a homology equivalence so comes from a homological
2-deformation. However the intermediate stages in a homological 2-deformation will usually not be
homotopy equivalences (the fundamental groups change) so further data is needed to encode the fact that
the composition is. The proof shows that the homological moves can be done
together, so there are maps
$$K@>>>K' @>>>L'@>>> L$$
whose composition is the given map, the first and last are 2-deformations, and $K'\to L'$ changes (all)
2-cell attaching maps by $Z\pi_1$ homology in the 1-skeleton. 

$Z\pi_1$ homology corresponds to surfaces whose $\pi_1$ maps trivially into both $K'$ and $L'$. Therefore
curves on the surface extend to maps of 2-disks into both $K'$ and $L'$. Conversely if there are
enough 2-disks on the surface to surger away the extra $\pi_1$ then this shows the attaching maps are
homotopic and the complexes are homotopy equivalent. Thus simple {\it homotopy\/} equivalence can be
characterized as 2-deformations, and simultaneous homologies of attaching maps together with appropriate
nullhomotopies of curves on the homology surfaces. An explicit description of this data is given in
\cite{Q}, where it is
 called an ``s-move.'' 

\head 7. Proof of the deformation theorem\endhead
The data is a decomposition $(N,\partial N)=(M,A)\cup(W,V)$, a map of a relative 2-complex $(K,A)\to
(N,A)$, and a simple equivalence of chain complexes over $Z\pi_1N$, $C_*^c(K,A)\to C_*^c(M,A)$. 

The first step is to use the 1-skeleton alignment lemma: change $(K,A)$ by relative 2-deformation and
$(M,A)$ by handle moves so the chain map becomes a basis-preserving isomorphism in degrees 0 and~1. Then
since the complexes are 2-dimensional and the map is an equivalence, it must be an isomorphism in degree
1 too. Since the chain map is simple and the other degrees are basis-preserving, this isomorphism is
simple. This means that after stabilization it is a product of a diagonal matrix with group elements on
the diagonal, and elementary matrices. The diagonal matrix can be realized by changing choice of lift to
the cover used as basis element. Elementary matrices can be realized by handle slides in $M$, or homotopy
of attaching maps in $K$. After doing this the chain map becomes a basis-preserving isomorphism in degree
2 as well.

Choose a homotopy equivalence of (relative) 1-skeletons $(K^1,A)\to (M^1,A)$ realizing the bijection of
basis elements in the cellular chain complexes. Use this to replace the 1-skeleton of $K$, and regard $K$
as obtained from $M^1$ by attaching 2-cells. The bijection of bases in
$C^c_*$ gives a correspondence between 2-cells of $(K,A)$ and 2-handles of $(M,A)$. Since the boundary
homomorphisms agree the attaching maps are homologous in the 1-skeleton.

Choose surfaces mapping into the 1-skeleton representing the homologies between attaching
maps of $K$-2-cells and $M$-2-handles. Choose symplectic basis curves for the fundamental groups of
these surfaces. These are collections of simple closed curves in the surfaces so that each curve
intersects exactly one other curve, in a single point. Since the homology takes place in the cover
corresponding to
$\pi_1N$ each of these curves is nullhomotopic in $N$. Choose nullhomotopies. These extend the surfaces to
maps of capped surfaces into $N$.
See \cite{FQ} for details about capped surfaces.  

The objective is to do ambient 2-deformations of $W$ to change the 2-handles of $M$ so their attaching
maps become homotopic to the 2-cells of $K$. The capped surfaces are used to see how to do this.

The surfaces map into the 1-skeleton. After homotopy the caps can be arranged to have subdisks that map to
parallels of the cores of the 2-handles, and the complement of these subdisks maps into $M^1$. A cap can
be split along an arc by joining the base surface to itself by a tube along the arc. The fragments give
caps for the new surface, with dual caps given by parallel copies of the dual to the original. Repeated
splitting gives a capped surface so that each cap has a single subdisk going over a 2-handle.

The next step is a subdivision trick in $W$ that reduces to the case where all caps go over (duals of)
2-handles in $W$, and each such 2-handle intersects a single cap. This is in preparation for
homotopically trivial handle moves in $W$ which will change the decomposition. In principle a lot could be
done without changing the decomposition, using handle moves in
$M$. Suppose there is a cap that goes over a handle $h_i$ in $M$, and this cap is on the homology surface
between an attaching map in $K$ and a handle $h_j$ in $M$, with $i\neq j$. The cap describes how to do
two handle moves of $h_j$ over $h_i$ to move past the cap and reduce the genus of the homology surface by
one. The remainder of the  surface shows this move does not change the chain complex data. However moving
$h_j$ disturbs all the caps passing over it,  adding new subdisks passing over
$h_i$.  Care would necessary to be sure we are making progress. Since we have to do the more disruptive
moves anyway, we use them for everything and dodge this point. 

We describe the subdivision trick. Suppose $h$ is a 2-handle of the decomposition of $(N,A)$ that does
not lie in $M$. In other words it is dual to a 2-handle in the structure on $(W,V)$, and is given by an
embedding $h\:(D^2\times D^2, S^1\times D^2)\to (W,\partial W-V)$. Suppose
$n$ caps pass over $h$. This means they intersect the image in disks $D^2\times\{p_i\}$ for $p_i$ points
in the other $D^2$ factor. Think of $D^2$ as a handlebody with $n$ 0-handles and $n-1$ 1-handles, for
instance with spine an interval subdivided to have  $n$ vertices. Identify (by isotopy) the $p_i$ with
the center points of the 0-handles. Taking the product of  $(D^2, S^1)$ with this structure gives a
subdivision of $h$ into $n$ 2-handles and $n-1$ 3-handles. This separates the caps so each passes over a
distinct 2-handle. Recall that $h$ is the dual of a handle in the structure on $(W,V)$. The dual of this
operation is to split $h^*$ into 1- and 2-handles, so the 2-handlebody condition is preserved.

The point here is that ``passing over''  is not symmetric under duality. If a disk passes over a 2-handle
then it is parallel to the core. In the dual the core becomes the transverse disk which does not ``go
over'' the dual. Deforming it down to the attaching region gives a picture of it ``going under'' the dual
in the sense that it cuts through the attaching region but does not enter the interior of the handle.
Getting confused about these dual pictures is a mistake that has been discovered by all the foremost
4-dimensional topologists.

There is a variation on the subdivision trick that applies when $h$ is a handle in $M$. Again suppose $n$
caps pass over it. Subdivide it into $n+1$ 2-handles and $n$ 3-handles so that each cap gets its own
2-handle and there is one left over. Shrinking $h$ down to the left-over handle changes $M$ by ambient
isotopy and puts the rest of the handles into $W$. There they are dual to $n$ each   1- and
2-handles in $(W,V)$.
$W$ has not been changed up to diffeomorphism, and in fact these new handles geometrically cancel.  

We now have caps going over distinct 2-handles lying in $W$. We further improve the surfaces after
introducing some notation.  Let
$S_j$ denote the surface with boundary the composition of loops $(\partial h_j)(\partial d_j)^{-1}$,
where $d$ is a  2-cell of
$K$ and $h$ is a 2-handle of $M$. Let $D_{j,*}$ denote the caps attached to $S_j$. We want:
\roster\item each surface is divided into pieces $S_j=S'_j\cup R_j$ that intersect in an arc;
\item $R_j$ is a disk containing $\partial d_j$; 
\item $S'_j$ contains $\partial h_j$ and the boundaries of all the cap disks;
\item the sets $S'_j\cup D_{j,*}$ are disjointly embedded; and
\item $S'_j$ is contained in $\partial M-A$.
\endroster

This is not hard because the $S_j$ have 1-dimensional spines and it suffices to embed the spines. In
detail first choose a tree in each
$S_j$ that joins 
$\partial h_j$ to the intersection point in each dual pair of cap boundaries. Let
$S'_j$ be a small neighborhood in $S_j$ of the union of the tree, $\partial h_j$, and the cap boundaries.
Denote by $R_j$ the closure of the complement; by construction this is a disk containing
$\partial d_j$. Denote by $T_j$ the tree union with an arc in each cap which connects the intersection
point in the cap boundaries to the subdisk mapping to the handle core. These start off in $M^1$ but since
they are 1-dimensional we can push them rel endpoints on $\partial h_j$ and handle cores to lie in
$\partial M-A$. Since this is 3-dimensional we  can approximate the map to be disjoint embeddings on
the $T_j$. We can can extend the embedding $T_j\subset \partial M-A$ to an  embedding of a neighborhood
of $T_j$ in  $S_j$, and thus an embedding $S'_j\subset \partial M-A$. 

The next step is to do handle moves. Choose one in each dual pair of caps, and use it to do surgery on
$S'_j$. Surgery consists of cutting out a neighborhood of the boundary of the cap and replacing it by two
parallel copies. We choose the parallel copies to be on the boundary of the handle whose core is the cap.
The net effect is to convert each $S_j$ into a disk embedded in the boundary of
$M^1\cup(\text{2-handles})$. Push the $h_j$ across these disks by isotopy. On the handlebody level this
induces handle slides of the $h_j$ over  cap 2-handles. Afterwards the attaching maps of the
handles are homotopic in the 1-skeleton to the $K$ 2-cell attaching maps, via the disks $R_j$. Therefore
the move changes $M$ to have 2-skeleton equal to $K$ (that is, equal to the complex currently denoted $K$,
and obtained from the original $K$ by 2-deformation).

The final step is to describe these handle moves from the dual point of view, and recognize them as
ambient 2-deformations of $W$. First recall that surgery on $S'_j$ using a cap yields two parallel copies
of the cap in the resulting disk. Changing the attaching map of $h_j$ by isotopy across the disk does two
handle additions of $h_j$ over the cap handle. The dual cap gives a homotopy between arcs recording the
homotopy classes of these additions. This homotopy lies in $N$ so shows the two additions cancel in the
chain complex with $Z\pi_1N$ coefficients. If the homotopy were in $M$ it would show the additions cancel
homotopically, and give an ambient 2-deformation of $M$ in $N$. However these homotopies go over handles
in $W$ and are not in $M$. Some can be deformed into $M$, but this is unusual.

The dual of a handle move of $h$ over $g$ is a move of $g^*$ over $h^*$. The dual of the procedure above
thus involves moves of 2-handles of $(W,V)$  over duals of handles in $(M,A)$. These moves still occur in
algebraically cancelling pairs, with a homotopy between the addition arcs provided by the dual cap.
Since the cap is in $W$ the homotopy also lies in $W$. Thus viewed from this side  the moves are
homotopically cancelling handle additions and therefore
ambient 2-deformations of
$W$. 

This completes the geometric part of the theorem. The ``Refinement'' asserts that when the geometric
output is converted back into algebra it matches up with the algebraic input. This is supposed to be
clear since at every step the geometric moves were modeled on the algebra. We omit details since they are
routine but long.

\bye